\documentclass[11pt]{article}
\usepackage{latexsym}
\usepackage{epsfig}
\usepackage{amsfonts,amssymb,amsmath,amsthm}
\usepackage{graphicx,multicol,color}
\usepackage{stmaryrd}
\usepackage{float}
\usepackage[utf8]{inputenc}
\usepackage[T1]{fontenc}
\usepackage{caption}

\usepackage{bold-extra}
\usepackage{color}                                            
\usepackage{array}                                            
\usepackage{longtable}                                        
\usepackage{calc}                                             
\usepackage{multirow}                                         
\usepackage{hhline}                                           
\usepackage{ifthen}    

\title{The Hausdorff dimension of the range of the Lévy multistable processes}
\author{R. Le Gu\'evel\\
\small{{\em Universit\'{e} Rennes 2 - Haute Bretagne, Equipe de Statistique Irmar, UMR CNRS 6625}}\\
\small{{\em Place du Recteur Henri Le Moal, CS 24307, 35043 RENNES Cedex, France}}\\
\small{{\em ronan.leguevel@univ-rennes2.fr}}}

\date{}

\def\bbbr{{\bf R}} 
\def\bbbn{{\bf N}} 

\newtheorem{theo}{Theorem}[section]

\newtheorem{lem}[theo]{Lemma}

\newcommand\li{\lim_{n \rightarrow +\infty}\limits}

\renewcommand\P{{\sf P}}

\newcommand{\one}{\ifmmode {\sf 1}\hspace{-.26em}{\sf
l}\hspace{-.35em}{\sf \_} \else ${\sf 1}\hspace{-.26em}{\sf
l}\hspace{-.35em}{\sf \_}$ \fi}

\renewcommand{\Box}{\mbox{\rule{1ex}{1ex}}}

\begin{document}

\maketitle

\begin{abstract}

We compute the Hausdorff dimension of the image $X(E)$ of a non random Borel set $E \subset [0,1]$, where $X$ is a L\'evy multistable process in $\bbbr.$
This extends the case where $X$ is a classical stable L\'evy process by letting the stability exponent $\alpha$ be a smooth function. Hence we are considering here non-homogeneous processes with increments which are not stationary and not necessarily independent. Contrary to the situation where the stability parameter is a constant, 
the dimension depends on the version of the multistable L\'evy motion when the process has an infinite first moment.
\end{abstract}

\section{Introduction}

For $(X_t)_t$ a stochastic process, we define the range of $X$ on a non random Borel set $E$ as the set $X(E) = \{ x : x=X_t \textrm{ for  some } t \textrm{ in } E \}$.
 We already know that for $X$ a typical L\'evy process, $X(E)$ is a random fractal set.
Many authors have been interested in producing the dimension properties of the sets $X(E)$. The computation of $\dim X(E)$ has been performed under various assumptions on $X$ and $E$, mainly if $X$ is a stable process, a subordinator or a general L\'evy process. For instance, see MacKean \cite{McKean},
Blumenthal and Getoor \cite{BGet}, Hawkes \cite{Haw}, Pruitt and Taylor \cite{PruTay}, Hendricks \cite{Hen} or Kahane \cite{Kah} for stable processes, Millar \cite{Millar}, Pruitt \cite{Pruitt} or Blumenthal and Getoor \cite{BGet2} for processes with stationary independent increments.
More recently, some results on operator stable processes or additive Lévy processes have been obtained for example in Becker-Kern, Meerschaert and Scheffler \cite{BMS}, Meerschaert and Xiao \cite{MX},
Khoshnevisan, Xiao and Zhong \cite{KhoXZ} or Khoshnevisan and Xiao \cite{KhoX}.
 Our aim in this article is to present the fractal nature of $X(E)$ through its Hausdorff dimension, with the assumption that $X$ belongs to the class of multistable Lévy processes, a natural extension of the stable Lévy processes.

The multistable processes have been introduced by Falconer and Lévy-Véhel in 2009 \cite{KFJLV}. 
Their distributions, their H\"olderian regularity or their multifractal properties have been studied for instance in
\cite{AA,LL,LL2,LL3,KFLL}. They provide useful models for all applications that deal with discontinuous processes where the intensity of jumps is non-constant. Most multistable processes are non-homogeneous in the sense that 
their increments are neither independent nor stationary. In this article, we consider only multistable L\'evy motions which are the simplest examples of multistable processes.

The paper is organised as follows. Section \ref{Sec2} contains the notations. In Section \ref{Sec3}, we present the main results on the computation of the Hausdorff dimension of the range.
Section \ref{Sec4} is dedicated to statement of useful technical lemmas on multistable processes. All the proofs are gathered in Section \ref{Sec5}.

\section{Notations}\label{Sec2}

We first summarise the basic notions about Hausdorff measures on the real line (see Falconer \cite{Falc} for more details).
For a subset $E$ of $[0,1]$, the diameter of $E$ is defined as $|E| = \sup \{|x-y| : x\in E, y \in E\}.$ Let 
$\beta$ be a non-negative number. For any $\delta >0$ we define
$$\mathcal{H}^{\beta}_{\delta}(E) = \inf \left\{ \sum_{i=1}^{+\infty}\limits |U_i|^{\beta} : \{ U_i\} \textrm{ is a } \delta-\textrm{cover of } E \right\}.$$
We call $\mathcal{H}^{\beta}(E) = \lim_{\delta \rightarrow 0}\limits \mathcal{H}^{\beta}_{\delta}(E)$ the $\beta$-dimensional Hausdorff measure of $E$, and the Hausdorff dimension of $E$ is defined as 
$$\dim(E) = \inf \left\{ \beta : \mathcal{H}^{\beta}(E) = 0 \right\} = \sup \left\{ \beta : \mathcal{H}^{\beta}(E) = \infty \right\}.$$
Throughout the paper, $\textrm{c}(E)$ stands for the convex hull of $E$, that is $\textrm{c}(E) = \{ tx +(1-t)y : t \in [0,1], x\in E, y \in E\}$.
$\mathring{E}$ will be the interior of $E$, and $\mathcal{P}$ will represent the set of partitions of $[0,1]$.
For $A \in \mathcal{P}$, we shall write $A=A^n$ if the number of intervals composing $A$ is $n$, and if $A^n=(A_i^n)_{i=1,...,n}$ is such that $[0,1] = \bigcup_{i=1}^{n}\limits A_i^n$ and $A_i^n \cap A_j^n = \varnothing$ for $i\neq j$, the mesh of $A^n$ is defined as
$|A^n| = \max_{i=1}^n\limits |A_i^n|$. Without loss of generality, $A_1^n$ is assumed to be the first set, that is for all $n \geq 1$, $0 \in A_1^n$.

The remainder of this section will be devoted to the definition of the multistable Lévy processes, using their Ferguson-Klass-LePage representation. 
For $\alpha \in (0,2]$, recall that the stochastic integral $I(f):= \int f(x) M(dx)$ of a real function $f$ with respect to $M$ exists if, for instance, $M$ is a symmetric $\alpha$-stable random measure on $\mathbb{R}$, with the Lebesgue measure as the control measure, and if $f$ is measurable and satisfies 
$\int_{\mathbb{R}} |f(x)|^{\alpha} dx < + \infty$ (see \cite{ST}). Write $S_{\alpha}(\sigma, \beta,\mu)$ for the $\alpha$-stable distribution with scale parameter $\sigma$, skewness $\beta$ and shift parameter $\mu$; many symmetric stable processes $\{ Y(t), t \in \mathbb{R}\}$ admit the stochastic integral representation
$$Y(t) = \int f_t(x) M(dx).$$
The marginal distribution of $Y$ is therefore $Y(t) \sim S_{\alpha}(\sigma_{f_t},0,0)$
where $\sigma_{f_t} = \left(\int_{\mathbb{R}} |f_t(x)|^{\alpha} dx\right)^{1/\alpha}.$

Among them, the symmetric standard $\alpha$-stable Lévy process on the interval $[0,1]$ may be defined as 
$$L_{\alpha}(t) := \int_{\mathbb{R}} \mathbf{1}_{[0,t]}(x) M(dx), \quad t\in [0,1].$$
Since $L_{\alpha}(t) \sim S_{\alpha}(t^{1/\alpha},0,0)$, the logarithm of the characteristic function of $L_{\alpha}(t)$ is given by $\log E[e^{ i \theta L_{\alpha}(t)}] = - t|\theta|^{\alpha}.$

We shall use another representation of the Lévy processes, based on series of random variables, in order to define its multistable versions.
We need for that the following sequences:

\begin{itemize}
	\item $(\Gamma_i)_{i \geq 1}$ a sequence of arrival times of a Poisson process with unit arrival rate,
	\item $(V_i)_{i \geq 1}$ a sequence of i.i.d. random variables with uniform distribution on $[0,1]$, independent of $(\Gamma_i)_{i \geq 1}$,
	\item $(\gamma_i)_{i \geq 1}$ a sequence of i.i.d. random variables with distribution $P(\gamma_i=1)=P(\gamma_i=-1)=1/2$, independent of $(\Gamma_i)_{i \geq 1}$ and $(V_i)_{i \geq 1}$.
\end{itemize}
Accordingly, the Lévy motion $\{ L_{\alpha}(t), t \in [0,1] \}$ admits the series representation:
$$
L_{\alpha}(t) = \sum_{i=1}^{\infty} \gamma_i C_{\alpha} \Gamma_i^{-1/\alpha} \mathbf{1}_{[0,t]}(V_i)
$$
where $C_{\alpha} = \left( \int_{0}^{\infty} x^{-\alpha} \sin x \ dx \right)^{-1/\alpha}.$ For more details about Ferguson-Klass-LePage representations, we refer the reader to \cite{FK, JR, ST}. 
It becomes clear that the stable Lévy motion is a càdlàg process, jumping at time $V_i$ with a jump of size $\Gamma_i^{-1/ \alpha}$, that is the stability index $\alpha$ may be seen as a parameter fitting the size of the jumps. 

The multistable processes are more flexible since they allow us to consider non homogeneous jumps processes with a non constant index of stability $\alpha$.
The size of the jumps will be governed by a function $\alpha(t)$ evolving with time. The first way to define such a process is to use the Ferguson-Klass-LePage representation of the stable processes, as in \cite{LL}, replacing $\alpha$ by a function $\alpha : [0,1] \rightarrow (0,2)$. From now on
we make the assumption that $\alpha$ is $\mathcal{C}^1$ function, ranging in $[\alpha_*,\alpha ^*]$, a subset of $(0,2)$.
The multistable L\'evy motion is the process
$$
X(t) = \sum_{i=1}^{\infty} \gamma_i C_{\alpha(t)} \Gamma_i^{-1/\alpha(t)} \mathbf{1}_{[0,t]}(V_i).
$$
Since we have replaced $\alpha$ by $\alpha(t)$, for each $t \in [0,1]$, $X(t)$ is a symmetric $\alpha(t)$-stable random variable $S_{\alpha(t)}(t^{1/\alpha(t)},0,0)$ and $\log E[e^{ i \theta X(t)}] = - t|\theta|^{\alpha(t)}.$

The second construction, due to Falconer and Liu \cite{KFLL}, comes from the definition of multistable random measures $M_{\alpha(x)}$ where we have replaced again $\alpha$ by a function $\alpha(t)$. They defined the stochastic integral of $f$ with respect to
a multistable random measure providing all its finite dimensional distributions. The multistable Lévy motion resulting from this definition is 

$$
Z(t) = \sum_{i=1}^{\infty} \gamma_i C_{\alpha(V_i)} \Gamma_i^{-1/\alpha(V_i)} \mathbf{1}_{[0,t]}(V_i),
$$
which satisfies $\log E[e^{ i \theta Z(t)}] = - \int_0^t|\theta|^{\alpha(x)} dx.$

We already know that the two processes $X$ and $Z$ are linked by the following formula (\cite{LGLVL}, Theorem 8):

\begin{equation}\label{semimart}
X(t) = Y(t) + Z(t),
\end{equation}
where $Y(t) = \int_0^t\limits \sum_{i=1}^{+\infty}\limits \gamma_i K_i(u) \mathbf{1}_{[0,u[}(V_i) du$ and $K_i(u) = \frac{d\left(C_{\alpha(s)} \Gamma_i^{-1/\alpha(s)}\right)}{ds}(u).$

Our results involve the following quantities: $\alpha_*(E)=\inf_{t \in E}\limits \alpha(t)$, $\alpha^*(E)=\sup_{t \in E}\limits \alpha(t)$,
$$d_*(E) = \max(1,\alpha_*(E)) \dim (E) \quad \textrm{and} \quad d^*(E) = \max(1,\alpha^*(E)) \dim (E).$$

Finally, in all the paper, for some parameter $\beta$, $K_{\beta}$ will mean a finite positive constant which depends only on $\beta$, and 
we will use the fact that there exists $K >0$ such that for all $u \in U$ and all $i \geq 1$, 
\begin{equation}\label{noyau}
|K_i(u)| \leq K(1+|\log \Gamma_i|)(\frac{1}{\Gamma_i^{1/\alpha_*(U)}} + \frac{1}{\Gamma_i^{1/\alpha^*(U)}}).
\end{equation}

\section{Main theorems}\label{Sec3}

\begin{theo}\label{Theo1}
Let $E$ be a subset of $[0,1]$. Almost surely,
$$\dim Z(E) \geq \min(1,\alpha_*(c(E)) \dim (E)) .$$

Suppose also that $\inf_{s \in E}\limits s >0$, $\sup_{(s,t) \in E^2}\limits \frac{|t-s|}{|\alpha(t)-\alpha(s)|} < +\infty$ and $\alpha^*(\textrm{c}(E)) - \alpha_*(\textrm{c}(E)) \leq \frac{\alpha_*^2}{2}$ . Almost surely, 
$$\dim X(E) \geq \min(1,d_*(E)) .$$
\end{theo}

\begin{theo}\label{Theo2}
Let $E$ be a subset of $[0,1]$. Almost surely, 
$$\dim Z(E) \leq \alpha^*(E) \dim (E)$$
and
$$\dim X(E) \leq d^*(E).$$
\end{theo}

\begin{theo}\label{Theo3}
Let $(A^n)_{n \in \mathbb{N}}$ be a sequence of partitions of $\mathcal{P}$ such that $\li |A^n|=0$, $E$ a subset of $[0,1]$.
Almost surely, 
$$\dim Z(E) = \min(1, \limsup_{n \rightarrow +\infty}\limits \max_{i=1}^{n}\limits \alpha^*(E \cap A_i^n) \dim(E \cap A_i^n) ).$$
\end{theo}

\begin{theo}\label{Theo4}
Let $(A^n)_{n \in \mathbb{N}}$ be a sequence of partitions of $\mathcal{P}$ such that $\li |A^n|=0$, $E$ a subset of $[0,1]$ such that $\inf_{s \in E}\limits s >0.$
Assume that $\exists n_0 \geq 1$ such that $\forall n \geq n_0$, $\forall i \in \llbracket 1,n \rrbracket$, $\sup_{(s,t) \in (E\cap \mathring{A_i^n})^2}\limits \frac{|t-s|}{|\alpha(t)-\alpha(s)|} < +\infty.$
Almost surely, 
$$\dim X(E) = \min(1, \limsup_{n \rightarrow +\infty}\limits \max_{i=1}^{n}\limits d^*(E \cap A_i^n)) = \min(1, \limsup_{n \rightarrow +\infty}\limits \max_{i=1}^{n}\limits d_*(E \cap A_i^n)).$$
\end{theo}

{\bf Remark:} Notice that if $\alpha^*(E) \leq 1$, almost surely, $\dim X(E) = \dim(E),$ which is not true if $\alpha$ is constant. Else, if $\alpha^*(E) > 1$, almost surely, $\dim X(E) = \dim Z(E).$

\section{Technical lemmas}\label{Sec4}

\begin{lem}\label{Lem1}
$\forall \beta \in (0,1)$, $\forall U \subset [0,1]$, $\exists K_{U,\beta} >0$ such that $\forall (s,t) \in U^2$,
\begin{eqnarray*}
\textrm{     } E[|X(t)-X(s) |^{-\beta} ] & \leq & K_{U,\beta} |t-s|^{-\frac{\beta}{\alpha_*(U)}}\\
\textrm{ and } E[|Z(t)-Z(s) |^{-\beta} ]& \leq & K_{U,\beta} |t-s|^{-\frac{\beta}{\alpha_*(c(U))}}.
\end{eqnarray*}
If we assume also that $\inf_{s \in U}\limits s >0$, $\sup_{(s,t) \in  U^2}\limits \frac{|t-s|}{|\alpha(t)-\alpha(s)|} < +\infty$ and $\alpha^*(\textrm{c}(U)) - \alpha_*(\textrm{c}(U)) \leq \frac{\alpha_*^2}{2}$, then $\exists K_{U,\beta} >0$ such that $\forall (s,t) \in U^2$,
$$E[|X(t)-X(s) |^{-\beta} ]\leq K_{U,\beta} |t-s|^{-\beta}.$$
\end{lem}

\begin{lem}\label{Lem2}
Let $(I_j)_{j}=([a_j,b_j])_{j}$ be a collection of closed intervals of $[0,1]$ and $p \in (0, \inf_{j}\limits \alpha_*(I_j))$. For all $\varepsilon >0$, $\exists K_{p,\varepsilon} >0$ such that
$\forall j$,
\begin{eqnarray*}
\textrm{     } E[\sup_{(s,t) \in I_j^2}|Z(t)-Z(s) |^{p} ]& \leq & K_{p,\varepsilon} |I_j|^{\frac{p}{\max(1,\sup_{j} \alpha^*(I_j))+\varepsilon}}\\
\textrm{ and } E[\sup_{(s,t) \in I_j^2}|X(t)-X(s) |^{p} ]& \leq & K_{p,\varepsilon} |I_j|^{\frac{p}{\max(1,\sup_{j} \alpha^*(I_j))+\varepsilon}}.
\end{eqnarray*}
\end{lem}

\begin{lem}\label{Lem3}
Let $(A^n)_{n \in \bbbn} \in \mathcal{P}^{\bbbn}$ be a sequence of partitions of $[0,1]$ such that $\li |A^n| = 0$. Then, for all subsets $E$ of $[0,1]$, 
\begin{eqnarray}\label{equalimsupd}
\limsup_{n \rightarrow +\infty}\limits \max_{i=1}^{n}\limits d^*(E \cap A_i^n) & =& \limsup_{n \rightarrow +\infty}\limits \max_{i=1}^{n}\limits d_*(E \cap A_i^n)
\end{eqnarray}
and
\begin{eqnarray*}
\limsup_{n \rightarrow +\infty}\limits \max_{i=1}^{n}\limits \alpha^*(E \cap A_i^n) \dim(E \cap A_i^n) & =&\limsup_{n \rightarrow +\infty}\limits \max_{i=1}^{n}\limits \alpha_*(E \cap A_i^n) \dim(E \cap A_i^n)\\
& = & \limsup_{n \rightarrow +\infty}\limits \max_{i=1}^{n}\limits \alpha^*(c(E \cap A_i^n)) \dim(E \cap A_i^n)\\
& = & \limsup_{n \rightarrow +\infty}\limits \max_{i=1}^{n}\limits \alpha_*(c(E \cap A_i^n)) \dim(E \cap A_i^n).
\end{eqnarray*}
Furthermore, all these equalities also occur with $\liminf_{n \rightarrow +\infty}\limits$ instead of $\limsup_{n \rightarrow +\infty}\limits$.

\end{lem}

\section{Proofs}\label{Sec5}

\noindent \textbf{\textsc{Proof of theorem \ref{Theo1}}}

Let $\beta < \min(1,\alpha_*(c(E))\dim(E)).$ Since $\frac{\beta}{\alpha_*(c(E))} < \dim(E)$, $\mathcal{H}^{\frac{\beta}{\alpha_*(c(E))}}(E) = +\infty.$
According to Davies theorem \cite{Dav}, $\exists F \subset E$, $F$ closed set such that $\mathcal{H}^{\frac{\beta}{\alpha_*(c(E))}}(F) >0.$ Then $C_{\frac{\beta}{\alpha_*(c(E))}} (F) >0$ by Frostman's theorem.
Let $p_m$ a probability measure concentrated on $F$ s.t. $$\int_F \int_F |x-y|^{-\frac{\beta}{\alpha_*(c(E))}} p_m(dx) p_m(dy) <+\infty.$$

With Lemma $\ref{Lem1}$,
\begin{eqnarray*}
E \left[ \int_F \int_F |Z(x)-Z(y)|^{-\beta} p_m(dx) p_m(dy) \right] &\leq& K_{\beta,F} \int_F \int_F |x-y|^{-\frac{\beta}{\alpha_*(c(F))}} p_m(dx) p_m(dy)\\
&\leq& K_{\beta,F} \int_F \int_F |x-y|^{-\frac{\beta}{\alpha_*(c(E))}} p_m(dx) p_m(dy)\\
& < & +\infty.
\end{eqnarray*}
So $\P \left( \mathcal{H}^{\beta}(Z(F)) >0 \right) =1$, $\P \left( \mathcal{H}^{\beta}(Z(E)) >0 \right) =1,$ and $\dim(Z(E)) \geq \beta$.

Assume now that $\inf_{s \in E}\limits s >0$, $\sup_{(s,t) \in E^2}\limits \frac{|t-s|}{|\alpha(t)-\alpha(s)|} < +\infty$ and $\alpha^*(\textrm{c}(E)) - \alpha_*(\textrm{c}(E)) \leq \frac{\alpha_*^2}{2}.$
The proof for the process $X$ is similar to the previous one. Consider $\beta < \min(1,d_*(E))$ and $\gamma_*(E) = \max(1,\alpha_*(E)).$ We obtain $\dim(X(E)) \geq \beta$ replacing $\alpha_*(c(E))\dim(E)$ by $d_*(E)$ in the previous calculus and $ \alpha_*(c(E))$ by $\gamma_*(E)$ $\Box$
%
%

\noindent \textbf{\textsc{Proof of theorem \ref{Theo2}}}

For a partition $(A_k)_{k=1,...N}$,
$\dim X(E) = \max_{k=1}^N\limits \dim X(E \cap A_k)$ therefore it is enough to show that for all $k$,
$$\dim X(E \cap A_k) \leq \max(1, \alpha^*(E \cap A_k)) \dim(E \cap A_k) \left( \leq \max(1,\alpha^*(E) )\dim E \right)$$
and
$$\dim Z(E \cap A_k) \leq \alpha^*(E \cap A_k) \dim(E \cap A_k) \left( \leq \alpha^*(E) \dim E \right).$$
Thus we may suppose that $|\alpha^*(E) - \alpha_*(E)| \leq \varepsilon$ for $\varepsilon >0$ as small as we want.

\vspace{0.2cm}
\noindent {\it Suppose first that $\dim(E) < 1$.}
 \vspace{0.2cm}
 
 Let $\beta \in (\dim(E), 1)$ and $n_0 \in \mathbb{N}$. For each $n \geq n_0$, let $\{ I_{i n}, i\geq 1 \}$ be a cover of $E$ by closed intervals such that $\li \sum_{i=1}^{+\infty}\limits |I_{in}|^{\beta}=0.$
This can be done since $\mathcal{H}^{\beta}(E)=0$. Suppose also that $\varepsilon$ is small enough to have $\beta < \frac{\inf_{i,n \geq n_0}\limits \alpha_*(I_{in})}{\inf_{i,n \geq n_0}\limits \alpha_*(I_{in}) + 2\varepsilon} <1$,
and that $\sup_{i,n\geq n_0}\limits \alpha^*(I_{in}) < \inf_{i,n \geq n_0}\limits \alpha_*(I_{in}) +\varepsilon.$ We shall denote $c=\inf_{i,n \geq n_0}\limits \alpha_*(I_{in})$ and $d=\sup_{i,n \geq n_0}\limits \alpha^*(I_{in})$.
Notice that for all $i$,$n$, $\beta (d+\varepsilon) <  \beta (c+2\varepsilon) <c.$

Now for each $n\geq n_0$, $\{ X(I_{i n}), i\geq 1 \}$ is a cover of $X(E)$, and $\{ Z(I_{i n}), i\geq 1 \}$ a cover of $Z(E)$. We consider two cases to finish the proof when $\dim(E) <1.$

\noindent {\it $(i)$: Case $\alpha^*(E) \geq 1$.}

We apply Lemma \ref{Lem2} to obtain
$$E\left[\sum_{i=1}^{+\infty}\limits |X(I_{in})|^{\beta (d +\varepsilon)} \right] \leq K_{d,\beta,\varepsilon} \sum_{i=1}^{+\infty}\limits |I_{in}|^{\beta}$$
and
 \begin{equation}\label{MajespZ}
  E\left[\sum_{i=1}^{+\infty}\limits |Z(I_{in})|^{\beta (d +\varepsilon)} \right] \leq K_{d,\beta,\varepsilon} \sum_{i=1}^{+\infty}\limits |I_{in}|^{\beta}.
 \end{equation}
Then for a subsequence of $n$'s approaching $\infty$, almost surely, $\li \sum_{i=1}^{+\infty}\limits |X(I_{in})|^{\beta (d +\varepsilon)} = 0,$ and $\dim X(E) \leq \beta (\sup_{i,n \geq n_0}\limits \alpha^*(I_{in})+\varepsilon).$ 
Letting $\varepsilon$ tend to $0$, then letting $n_0$ tend to infinity one finally obtains $\dim X(E) \leq \beta \alpha^*(E).$ Since $\beta$ was arbitrary, $\dim X(E) \leq \alpha^*(E) \dim(E).$ Equation (\ref{MajespZ}) leads also to 
$\dim Z(E) \leq \alpha^*(E) \dim(E)$ for the same reasons.

\noindent {\it $(ii)$: Case $\alpha^*(E) <1$.}

Suppose that $\forall i$, $\forall n \geq n_0$, $\alpha^*(I_{in}) + \varepsilon <1.$ With equations $(\ref{semimart})$ and $(\ref{noyau})$,
$$|X(I_{in})| = \sup_{(s,t) \in I_{in}^2}\limits |X(t)-X(s)| \leq \sup_{(s,t) \in I_{in}^2}\limits \int_{s}^{t}\limits \sum_{j=1}^{+\infty}\limits K(1+|\log \Gamma_j|)(\frac{1}{\Gamma_j^{1/c}} + \frac{1}{\Gamma_j^{1/d}}) ds +  |Z(I_{in})|,$$

so
\begin{equation}\label{majdiametre}
\sum_{i=1}^{+\infty}\limits |X(I_{in})|^{\beta}  \leq \left(\sum_{j=1}^{+\infty}\limits K(1+|\log \Gamma_j|)(\frac{1}{\Gamma_j^{1/c}} + \frac{1}{\Gamma_j^{1/d}})\right)^{\beta} \sum_{i=1}^{+\infty}\limits |I_{in}|^{\beta} +  \sum_{i=1}^{+\infty}\limits |Z(I_{in})|^{\beta}.
\end{equation}
Since $\sum_{j=1}^{+\infty}\limits K(1+|\log \Gamma_j|)(\frac{1}{\Gamma_j^{1/c}} + \frac{1}{\Gamma_j^{1/d}})<+\infty$ and $\li \sum_{i=1}^{+\infty}\limits |I_{in}|^{\beta}=0$, almost surely, 
$$\li \left(\sum_{j=1}^{+\infty}\limits K(1+|\log \Gamma_j|)(\frac{1}{\Gamma_j^{1/c}} + \frac{1}{\Gamma_j^{1/d}})\right)^{\beta} \sum_{i=1}^{+\infty}\limits |I_{in}|^{\beta} =0.$$
Let us show that $\li \sum_{i=1}^{+\infty}\limits |Z(I_{in})|^{\beta} =0$ where the convergence is in probability.
$$\left| Z(t) - Z(s) \right| \leq K \sum_{j=1}^{+\infty}\limits (\frac{1}{\Gamma_j^{1/c}} + \frac{1}{\Gamma_j^{1/d}}) \mathbf{1}_{[s,t]}(V_j).$$

Let $D_{\alpha}(t) = \sum_{j=1}^{+\infty}\limits \frac{1}{\Gamma_j^{1/\alpha}} \mathbf{1}_{[0,t]}(V_j)$ so that 
\begin{equation}\label{inegalsubordinator}
|Z(I_{in})|^{\beta} \leq K_{\beta} |D_{c}(I_{in})|^{\beta} +  K_{\beta} |D_{d}(I_{in})|^{\beta}.
\end{equation}
$D_{d}$ is a stable-subordinator so $|D_{d}(I_{in})|^{\beta}$ is distributed as $|I_{in}|^{\beta/d} |D_{d}(1)|^{\beta}$.
Since $\frac{\beta}{d}> \beta$, $ \sum_{i=1}^{+\infty}\limits |D_{d}(I_{in})|^{\beta}$ tends to $0$ in probability. For the same reasons, $\li \sum_{i=1}^{+\infty}\limits |D_{c}(I_{in})|^{\beta} =0$ in probability, which entails with $(\ref{majdiametre})$ that $\li \sum_{i=1}^{+\infty}\limits |X(I_{in})|^{\beta} \stackrel{P}{=} 0.$
Then for a subsequence of $n$'s approaching $\infty$, almost surely, $\li \sum_{i=1}^{+\infty}\limits |X(I_{in})|^{\beta} = 0,$ and $\dim X(E) \leq \beta.$ Since $\beta$ was arbitrary,
$\dim X(E) \leq \dim(E).$ 

Replacing $\beta$ by $\beta d$ in the equation $(\ref{inegalsubordinator})$, we obtain $\li \sum_{i=1}^{+\infty}\limits |Z(I_{in})|^{\beta d} \stackrel{P}{=} 0,$ and $\dim Z(E) \leq \beta d.$
Letting $n_0$ tend to infinity one finally obtains $\dim Z(E) \leq \beta \alpha^*(E).$

\vspace{0.2cm}
\noindent {\it Suppose now that $\dim(E) = 1$.}
\vspace{0.2cm}

The result is obvious for the process $X$, and for the process $Z$ if $\alpha^*(E) \geq 1$ so we consider only the case $\alpha^*(E) <1$. As previously, the result is a consequence of the equation (\ref{inegalsubordinator}).
Let $\beta >1$, $n_0 \in \mathbb{N}$, $n \geq n_0$ and $\{ I_{i n}, i\geq 1 \}$ be a cover of $E$ by closed intervals such that $\li \sum_{i=1}^{+\infty}\limits |I_{in}|^{\beta}=0.$
Suppose also that $d = \sup_{i,n \geq n_0}\limits \alpha^*(I_{in}) <1.$ Equation $(\ref{inegalsubordinator})$ and its consequences are still available:
$$|Z(I_{in})|^{\beta d} \leq K_{\beta,d} |D_{c}(I_{in})|^{\beta d} +  K_{\beta,d} |D_{d}(I_{in})|^{\beta d}$$
leads to 
$\li \sum_{i=1}^{+\infty}\limits |Z(I_{in})|^{\beta d} \stackrel{P}{=} 0,$ $\dim Z(E) \leq \beta d,$ and $\dim Z(E) \leq \alpha^*(E) \Box$

\noindent \textbf{\textsc{Proof of Theorem \ref{Theo3} and Theorem \ref{Theo4}}}

Let us prove Theorem \ref{Theo4} first.
Suppose that $0 \in A_1^n$ for all $n \geq 1$. Since $\inf_{s \in E}\limits s >0$, for $n$ large enough, $E \cap A_1^n = \varnothing$.
We use Theorem $\ref{Theo2}$ to obtain $$\dim X(E) = \max_{i=1}^n\limits \dim X(E\cap A_i^n) \leq \max_{i=1}^n\limits d^*(E \cap A_i^n)$$
and
\begin{equation}\label{Firstineg}
\dim X(E) \leq \liminf_{n \rightarrow +\infty}\limits \max_{i=1}^n\limits d^*(E \cap A_i^n) \leq \limsup_{n \rightarrow +\infty}\limits \max_{i=1}^n\limits d^*(E \cap A_i^n).
\end{equation}
Let us show that $\dim X(E) \geq \min(1,\limsup_{n \rightarrow +\infty}\limits \max_{i=1}^n\limits d^*(E \cap A_i^n)).$  Theorem $\ref{Theo1}$ gives
\begin{equation}\label{Secondineg}
\dim X(E) = \max_{i=2}^n\limits \dim X(E\cap A_i^n) \geq \max_{i=2}^n\limits \min(1,d_*(E \cap A_i^n)).
\end{equation}
Then we consider three cases.

{\it $(i)$: Case $\limsup_{n \rightarrow +\infty}\limits \max_{i=1}^n\limits d^*(E \cap A_i^n) <1.$}

With the two inequalities $(\ref{Firstineg})$ and $(\ref{Secondineg})$, $\max_{i=1}^n\limits \min(1,d_*(E \cap A_i^n)) <1$ so for all $n \geq 1$ and all $i=1,...,n$,
$d_*(E \cap A_i^n) <1$ and $\max_{i=1}^n\limits \min(1,d_*(E \cap A_i^n)) = \max_{i=1}^n\limits d_*(E \cap A_i^n),$ i.e.
$$\dim X(E) \geq \max_{i=1}^n\limits d_*(E \cap A_i^n) .$$
Finally, $\dim X(E) \geq \limsup_{n \rightarrow +\infty}\limits \max_{i=1}^n\limits d_*(E \cap A_i^n)$ and the result comes from Lemma $\ref{Lem3}$.

{\it $(ii)$: Case $\limsup_{n \rightarrow +\infty}\limits \max_{i=1}^n\limits d^*(E \cap A_i^n) =1.$}

If for all $n \geq 1$ and all $i=1,...,n$,  $d_*(E \cap A_i^n) <1$, we obtain as previously $$\dim X(E) \geq \limsup_{n \rightarrow +\infty}\limits \max_{i=1}^n\limits d_*(E \cap A_i^n) =1.$$
Otherwise, there exists $n_0 \in \bbbn$ and $i_0 \in  \llbracket 1, n_0 \rrbracket$ such that $d_*(E \cap A_{i_0}^{n_0}) \geq 1.$  Then 
\begin{eqnarray*}
\dim X(E) & \geq & \dim X(E \cap A_{i_0}^{n_0})\\
& \geq & \min(1,d_*(E \cap A_{i_0}^{n_0}))\\
& = & 1.
\end{eqnarray*}

{\it $(iii)$: Case $\limsup_{n \rightarrow +\infty}\limits \max_{i=1}^n\limits d^*(E \cap A_i^n) >1.$}

With Lemma $\ref{Lem3}$, $\limsup_{n \rightarrow +\infty}\limits \max_{i=1}^n\limits d_*(E \cap A_i^n) >1$ so there exists $n_0 \in \bbbn$ and $i_0 \in  \llbracket 1, n_0 \rrbracket$ such that $d_*(E \cap A_{i_0}^{n_0}) \geq 1.$
As previously stated, $\dim X(E) \geq 1.$

In order to get Theorem \ref{Theo3}, replace $X$ by $Z$ and $d(E \cap A_i^n)$ by $\alpha(E\cap A_i^n) \dim(E \cap A_i^n)$ in the proof of Theorem \ref{Theo4} $\Box$

{\bf Remark:} Notice that if $\dim X(E) = \limsup_{n \rightarrow +\infty}\limits \max_{i=1}^n\limits d^*(E \cap A_i^n) <1$, then $\li \max_{i=1}^n\limits d^*(E \cap A_i^n)$ exists and is equal to $\dim X(E)$: indeed the inequality $(\ref{Firstineg})$
 becomes 
 $$\dim X(E) \leq \liminf_{n \rightarrow +\infty}\limits \max_{i=1}^n\limits d^*(E \cap A_i^n) \leq \limsup_{n \rightarrow +\infty}\limits \max_{i=1}^n\limits d^*(E \cap A_i^n)=dim X(E).$$
Lemma $\ref{Lem3}$ gives also in that case $\dim X(E) = \li \max_{i=1}^n\limits d_*(E \cap A_i^n).$

\noindent \textbf{\textsc{Proof of Lemma \ref{Lem1}}}


By Proposition 6.1 of \cite{LL}, the logarithm of the characteristic function of $X(t)-X(s)$ satisfies for $s \leq t$:
$$
\log \phi_{X(t)-X(s)}( \theta) = -2 s \int_0^{\infty}\limits \sin^2( \frac{\theta}{2} [\frac{C_{\alpha(t)}}{y^{1/\alpha(t)}} - \frac{C_{\alpha(s)}}{y^{1/\alpha(s)}} ]) dy - (t-s)|\theta|^{\alpha(t)},
$$
and by Proposition 2 of \cite{LGLVL}, 
$$
\log \phi_{Z(t)-Z(s)}( \theta) = - \int_s^t\limits |\theta|^{\alpha(u)} du.
$$
Accordingly for $|\theta| \geq 1$ and $(s,t) \in U^2$, 
\begin{equation}\label{caracteristicX}
\phi_{X(t)-X(s)}( \theta) \leq e^{-|t-s| |\theta|^{\alpha_*(U)}},
\end{equation}
\begin{equation}\label{caracteristicZ}
\phi_{Z(t)-Z(s)}( \theta) \leq e^{-|t-s| |\theta|^{\alpha_*([s,t])}} \leq e^{-|t-s| |\theta|^{\alpha_*(c(U))}},
\end{equation}
and for all $\theta$,
$$\phi_{X(t)-X(s)}( \theta) \leq e^{-2 \min(s,t) \int_0^{\infty}\limits \sin^2( \frac{\theta}{2} [\frac{C_{\alpha(t)}}{y^{1/\alpha(t)}} - \frac{C_{\alpha(s)}}{y^{1/\alpha(s)}} ]) dy}.$$
We obtain then for $(s,t) \in U^2$, using the Parseval's formula:
 \begin{eqnarray*}
  |t-s|^{\frac{\beta}{\alpha_*(U)}} E[|X(t)-X(s)|^{-\beta} ] & = & \int_0^{\infty}\limits \P (|X(t)-X(s)| \leq \frac{|t-s|^{\frac{1}{\alpha_*(U)}}}{x^{1/\beta}} ) dx\\
  & \leq & 1+ \frac{1}{\pi} \int_1^{\infty}\limits \int_{\bbbr}\limits \frac{\sin(\frac{\xi |t-s|^{\frac{1}{\alpha_*(U)}}}{x^{1/\beta}} )}{\xi} \phi_{X(t)-X(s)} (\xi) d\xi dx\\
  & = & 1+ \frac{1}{\pi} \int_1^{\infty}\limits \int_{\bbbr}\limits \frac{\sin(\frac{\theta}{x^{1/\beta}})}{\theta} \phi_{X(t)-X(s)} (\frac{\theta}{|t-s|^{\frac{1}{\alpha_*(U)}}}) d\theta dx\\
  & \leq & 1+ \frac{1}{\pi} \left( \int_1^{\infty}\limits \frac{dx}{x^{1/\beta}}\right)  \int_{\bbbr}\limits \phi_{X(t)-X(s)} (\frac{\theta}{|t-s|^{\frac{1}{\alpha_*(U)}}}) d\theta\\
  & \leq & 1+ \frac{1}{\pi} \frac{\beta}{1 - \beta} (2 |t-s|^{\frac{1}{\alpha_*(U)}} + 2 \int_{|t-s|^{\frac{1}{\alpha_*(U)}}}^{\infty}\limits e^{-|\theta|^{\alpha_*(U)}} d\theta )\\
  & \leq & 1+ \frac{1}{\pi} \frac{\beta}{1 - \beta}(2 + 2 \int_{0}^{\infty}\limits e^{-|\theta|^{\alpha_*(U)}} d\theta ).
 \end{eqnarray*}
 Using the same inequalities and $(\ref{caracteristicZ})$ instead of $(\ref{caracteristicX})$, we obtain also $$|t-s|^{\frac{\beta}{\alpha_*(c(U))}} E[|Z(t)-Z(s)|^{-\beta} ] \leq  1+ \frac{1}{\pi} \frac{\beta}{1 - \beta}(2 + 2 \int_{0}^{\infty}\limits e^{-|\theta|^{\alpha_*(c(U))}} d\theta ).$$

Assume now that $\inf_{s \in U}\limits s >0$, $\sup_{(s,t) \in  U^2}\limits \frac{|t-s|}{|\alpha(t)-\alpha(s)|} < +\infty$ and $\alpha^*(\textrm{c}(U)) - \alpha_*(\textrm{c}(U)) \leq \frac{\alpha_*^2}{2}$. 
Notice that $C_{\alpha(t)} = h \circ \alpha(t) $ where $h(v)= \left( \int_{0}^{\infty} x^{-v} \sin x \ dx \right)^{-1/v}$ is a continuously differentiable function on $[\alpha_*,\alpha^*]$.
Property 1.2.15 of \cite{ST} gives an explicit formula of $h$.  Then there exists $\omega_y \in [\alpha(t),\alpha(s)]$ (or $[\alpha(s), \alpha(t)]$) such that 
$$\frac{C_{\alpha(t)}}{y^{1/\alpha(t)}} - \frac{C_{\alpha(s)}}{y^{1/\alpha(s)}}  = (\alpha(t) - \alpha(s) ) \left[ \frac{h'(\omega_y) + h(\omega_y) \frac{\log(y)}{\omega_y^2} }{y^{1/ \omega_y}}\right].$$
Now the previous calculus gives 
\begin{eqnarray*}
 |t-s|^{\beta} E[|X(t)-X(s)|^{-\beta} ] & \leq & 1+ \frac{1}{\pi} \left( \int_1^{\infty}\limits \frac{dx}{x^{1/\beta}}\right)  \int_{\bbbr}\limits \phi_{X(t)-X(s)} (\frac{\theta}{|t-s|}) d\theta\\
 & \leq & 1+ K_{\beta} \int_{\bbbr}\limits e^{-2 \nu \int_0^{\infty}\limits \sin^2( \frac{\theta}{2} (\frac{\alpha(t) - \alpha(s)}{|t-s|} ) [ \frac{\omega_y^2 h'(\omega_y) + h(\omega_y) \log(y) }{\omega_y^2 y^{1/ \omega_y}} ]) dy}d\theta
\end{eqnarray*}
where $\nu = \inf_{s \in U}\limits s >0.$ Changing the variable $|\theta|$ according to the formula $\xi = \theta \frac{\alpha(t) - \alpha(s)}{|t-s|}$ leads to
$$|t-s|^{\beta} E[|X(t)-X(s)|^{-\beta} ]  \leq 1+ K_{\beta} \sup_{(s,t) \in  U^2}\limits \frac{|t-s|}{|\alpha(t) - \alpha(s)|}  \int_{\bbbr}\limits  e^{-2 \nu \int_0^{\infty}\limits \sin^2( \frac{\xi}{2}  [ \frac{\omega_y^2 h'(\omega_y) + h(\omega_y) \log(y) }{\omega_y^2 y^{1/ \omega_y}} ]) dy}d\xi.$$

Let $\varepsilon \in (0, \frac{\alpha_*^2}{4}).$ Using the fact that for $|x|$ small enough, $\sin^2(x) \geq \frac{1}{2} x^2,$ and the inequality $\inf_{x \in [\alpha_*,\alpha^*]}\limits |h(x)| >0,$ we may choose $K_1 >1$ and $K_2 >1$ such that for all $|\xi| \geq 1$, 
$$y \geq K_1 |\xi|^{\frac{\alpha^*(\textrm{c}(U))}{1-\varepsilon}} \Rightarrow \sin^2( \frac{\xi}{2}  [ \frac{\omega_y^2 h'(\omega_y) + h(\omega_y) \log(y) }{\omega_y^2 y^{1/ \omega_y}} ]) \geq K_2 |\xi|^{2} y^{-\frac{2}{\alpha_*(\textrm{c}(U))}}.$$
 Now 
 \begin{eqnarray*}
  \int_0^{\infty}\limits \sin^2( \frac{\xi}{2}  [ \frac{\omega_y^2 h'(\omega_y) + h(\omega_y) \log(y) }{\omega_y^2 y^{1/ \omega_y}} ]) dy & \geq & K_2 |\xi|^{2} \int_{y \geq K_1 |\xi|^{\frac{\alpha^*(\textrm{c}(U))}{1-\varepsilon}} }\limits y^{-\frac{2}{\alpha_*(\textrm{c}(U))}} dy\\
  & \geq & K |\xi|^{2+(1-\frac{2}{\alpha_*(\textrm{c}(U))} ) \frac{\alpha^*(\textrm{c}(U))}{1-\varepsilon}}.
 \end{eqnarray*}
 
 Since $ \frac{\alpha^*(\textrm{c}(U))}{\alpha_*(\textrm{c}(U))} \leq 1+ \frac{\alpha_*}{2}$, $2+(\frac{\alpha_*(\textrm{c}(U))-2}{\alpha_*(\textrm{c}(U))} ) \frac{\alpha^*(\textrm{c}(U))}{1-\varepsilon} >2+\frac{(\alpha_*-2)(2+\alpha_*)}{2(1-\varepsilon)} = \frac{\alpha_*^2-4\varepsilon}{2(1-\varepsilon)}.$ 
 Then for $|\xi| \geq 1$, $|\xi|^{2+(1-\frac{2}{\alpha_*(\textrm{c}(U))} ) \frac{\alpha^*(\textrm{c}(U))}{1-\varepsilon}} \geq |\xi|^{\frac{\alpha_*^2-4\varepsilon}{2(1-\varepsilon)}}$ and $\int_{\bbbr}\limits  e^{-2 \nu \int_0^{\infty}\limits \sin^2( \frac{\xi}{2}  [ \frac{\omega_y^2 h'(\omega_y) + h(\omega_y) \log(y) }{\omega_y^2 y^{1/ \omega_y}}]) dy}d\xi \leq K_{\varepsilon, U} <+\infty \Box$

 \noindent \textbf{\textsc{Proof of Lemma \ref{Lem2}}}

 Let $p \in (0, \inf_{j}\limits \alpha_*(I_j))$, $\varepsilon >0$ and $n_0 \in \bbbn$ large enough to have $n_0 \alpha_* > 2.$ Let $p' \in ( \max(1,\sup_{j}\limits \alpha^*(I_j)), 2),$ $c=\inf_{j}\limits \alpha_*(I_j)$ and $d=\sup_{j}\limits \alpha^*(I_j).$
 Equation $(\ref{semimart})$ can be written
 $$X(t) = \int_{0}^{t}\limits W_1(u) du + \int_{0}^{t}\limits W_2(u) du + Z(t)$$
 with $W_1(u) = \sum_{i=1}^{n_0}\limits \gamma_i K_i(u) \mathbf{1}_{[0,u[}(V_i).$ Then there exists a constant $K >0$ such that
 $$\sup_{(s,t)\in I_j^2}\limits |X(t)-X(s)|^p \leq K \left( \int_{a_j}^{b_j}\limits |W_1(u)| du \right)^p + K \left( \int_{a_j}^{b_j}\limits |W_2(u)| du\right)^p + K \sup_{(s,t)\in I_j^2}\limits |Z(t)-Z(s)|^p.$$
 
 The end of the proof consists of showing the inequality for these three terms. For the first term,
 $\int_{a_j}^{b_j}\limits |W_1(u)| du \leq (b_j-a_j) \sum_{i=1}^{n_0}\limits \sup_{u \in I_j}\limits |K_i(u)|$ so inequality (\ref{noyau}) gives:
 \begin{eqnarray*}
  ( \int_{a_j}^{b_j}\limits |W_1(u)| du )^p & \leq & K_{n_0} |b_j - a_j|^p \sum_{i=1}^{n_0}\limits \sup_{u \in I_j}\limits |K_i(u)|^p\\
   & \leq & K_p |b_j - a_j|^p \sum_{i=1}^{n_0}\limits(1+|\log \Gamma_i|)^p (\frac{1}{\Gamma_i^{1/c}} + \frac{1}{\Gamma_i^{1/d}} )^p,
 \end{eqnarray*}
hence $E \left[( \int_{a_j}^{b_j}\limits |W_1(u)| du )^p  \right]  \leq K_{n_0,p} |I_j|^p.$
 For the second term, we obtain by H\"older and Jensen inequalities
 \begin{eqnarray*}
  E \left[ ( \int_{a_j}^{b_j}\limits |W_2(u)| du )^p \right] & \leq & E \left[ ( \int_{a_j}^{b_j}\limits |W_2(u)| du )^{p'} \right]^{\frac{p}{p'}}\\
  & \leq & |b_j-a_j|^p \left( \sup_{u \in I_j}\limits E[|W_2(u) |^{p'}] \right)^{\frac{p}{p'}}.
 \end{eqnarray*}

 Since $K_i(u) \mathbf{1}_{[0,u[}(V_i)$ is independent of $\gamma_i$, we obtain with Theorem 2 of \cite{BE} that for all $u \in I_j$, 
 $$E[|W_2(u) |^{p'}] \leq \sum_{i >n_0}\limits E[|K_i(u)|^{p'} ].$$
 Then inequality (\ref{noyau}) leads to $\sup_{j}\limits \left( \sup_{u \in I_j}\limits E[|W_2(u) |^{p'}] \right)^{\frac{p}{p'}} < +\infty.$
 
 
 Let us consider the process $Z$. Proposition $5$ of \cite{LGLVL} yields that $Z$ is a semi-martingale and can be decomposed into $A+M$ where $M$ is a martingale and 
 $$M(t) = \sum_{i=1, \Gamma_i \geq 1}^{\infty}\limits \gamma_i C_{\alpha(V_i)} \Gamma_i^{-1/\alpha(V_i)} \mathbf{1}_{[0,t]}(V_i).$$
 Let $N = \textrm{Card} \{ i \geq 1 | \Gamma_i <1\}$ and $K_i = C_{\alpha(V_i)} \Gamma_i^{-1/\alpha(V_i)}$. We will use the following inequality: if $V_i \in I_j$, $K_i \leq K (\frac{1}{ \Gamma_i^{1/c}}+\frac{1}{ \Gamma_i^{1/d}})$ for some constant $K$ and the fact that $N$ is distributed as a Poisson random variable with unit mean. 
 For all $(s,t) \in [a_j,b_j]^2$,
 \begin{eqnarray*}
 |A(t)-A(s)|^p & = & \sum_{n=0}^{+\infty}\limits |\sum_{i=1}^{n}\limits \gamma_i K_i \mathbf{1}_{[s,t]}(V_i) |^p \mathbf{1}_{N=n}\\
 & \leq & K \sum_{n=1}^{+\infty}\limits n ( \frac{1}{ \Gamma_1^{1/c}}+\frac{1}{ \Gamma_1^{1/d}})^p (\sum_{i=1}^{n}\limits \mathbf{1}_{[a_j,b_j]}(V_i) )\mathbf{1}_{N=n}.
 \end{eqnarray*}
Using the fact that $V_i$ is independent of $\Gamma_1$ and $N$, 
   \begin{eqnarray*}
    E\left[\sup_{(s,t) \in I_j^2}\limits  |A(t)-A(s)|^p\right] & \leq & K(b_j-a_j) \sum_{n=1}^{+\infty}\limits n^2 E\left[ ( \frac{1}{ \Gamma_1^{1/c}}+\frac{1}{ \Gamma_1^{1/d}})^p \mathbf{1}_{N=n}\right].
   \end{eqnarray*}
Since $\sum_{n=1}^{+\infty}\limits n^2 E\left[ ( \frac{1}{ \Gamma_1^{1/c}}+\frac{1}{ \Gamma_1^{1/d}})^p \mathbf{1}_{N=n}\right] < +\infty$, $     E\left[\sup_{(s,t) \in I_j^2}\limits  |A(t)-A(s)|^p\right] \leq K(b_j-a_j).$
The last step of the proof is to show the inequality for the martingale $M$. 
Let $p'=\max(1,d)+\varepsilon.$ We apply the H\"older inequality to get $E[\sup_{(s,t) \in I_j^2}\limits  |M(t)-M(s)|^{p}] \leq E[\sup_{(s,t) \in I_j^2}\limits  |M(t)-M(s)|^{p'}]^{\frac{p}{p'}}.$
By the Doob's martingale inequality, there exists $K_{p'} >0$ such that 
$$E\left[\sup_{(s,t) \in I_j^2}\limits  |M(t)-M(s)|^{p'}\right] \leq K_{p'} \sup_{(s,t) \in I_j^2}\limits E[|M(t)-M(s)|^{p'} ].$$
 Now for every $(s,t)\in I_j^2$, Theorem 2 of \cite{BE} leads again to
 \begin{eqnarray*}
  E[|M(t)-M(s)|^{p'} ] & \leq & \sum_{i\geq 1}\limits E[|K_i|^{p'} \mathbf{1}_{\Gamma_i \geq 1} \mathbf{1}_{[s,t]}(V_i) ]\\
  & \leq & K_{p'}|b_j-a_j| \sum_{i=1}^{+\infty}\limits E[(\frac{1}{ \Gamma_i^{1/c}}+\frac{1}{ \Gamma_i^{1/d}})^{p'}\mathbf{1}_{\Gamma_i \geq 1}].
 \end{eqnarray*}
Since $\sum_{i=1}^{+\infty}\limits E[(\frac{1}{ \Gamma_i^{1/c}}+\frac{1}{ \Gamma_i^{1/d}})^{p'}\mathbf{1}_{\Gamma_i \geq 1}] < +\infty $, $(\sup_{(s,t) \in I_j^2}\limits E[|M(t)-M(s)|^{p'} ])^{\frac{p}{p'}} \leq K_{p,\varepsilon} |b_j-a_j|^{\frac{p}{p'}}$ which is the result of the Lemma \Box

 \noindent \textbf{\textsc{Proof of Lemma \ref{Lem3}}}
 
 Notice that 
 $$\left| d_*(E \cap A_i^n) - d^*(E \cap A_i^n) \right| \leq \left| (\alpha_*(E \cap A_i^n) - \alpha^*(E \cap A_i^n) \right| \dim( E \cap A_i^n).$$
 $\alpha$ is a $\mathcal{C}^1$ function so there exists $K >0$ such that the following inequalities hold:
 \begin{eqnarray}
 \label{majsupinf}  \left| d_*(E \cap A_i^n) - d^*(E \cap A_i^n) \right| & \leq & K \left| A^n \right|,\\
 \label{majsupinf2}  \left| \alpha_*(E \cap A_i^n) - \alpha^*(E \cap A_i^n) \right| \dim( E \cap A_i^n) &\leq & K \left| A^n \right|,\\
 \label{majsupinf3}   \left| \alpha_*(c(E \cap A_i^n)) - \alpha^*(c(E \cap A_i^n)) \right| \dim( E \cap A_i^n) &\leq& K \left| A^n \right|,\\
 \label{majsupinf4}    \left| \alpha^*(c(E \cap A_i^n)) - \alpha^*(E \cap A_i^n) \right| \dim( E \cap A_i^n) &\leq& K \left| A^n \right|.
 \end{eqnarray}

 Then, in order to prove equality (\ref{equalimsupd}), we use the inequality (\ref{majsupinf}) to obtain $d^*(E \cap A_i^n) \leq K \left| A^n \right| + d_*(E \cap A_i^n)$. This implies that
 $$\limsup_{n \rightarrow +\infty}\limits \max_{i=1}^{n}\limits d^*(E \cap A_i^n) \leq \limsup_{n \rightarrow +\infty}\limits \max_{i=1}^{n}\limits d_*(E \cap A_i^n)$$
 and 
 $$\liminf_{n \rightarrow +\infty}\limits \max_{i=1}^{n}\limits d^*(E \cap A_i^n) \leq \liminf_{n \rightarrow +\infty}\limits \max_{i=1}^{n}\limits d_*(E \cap A_i^n) .$$
 Equality (\ref{equalimsupd}) comes from the fact that $d_* \leq d^*.$ To obtain the second result of Lemma \ref{Lem3}, we may replace $d$ by $\alpha$ using  $(\ref{majsupinf2})$, $(\ref{majsupinf3})$ and $(\ref{majsupinf4})$ instead of $(\ref{majsupinf})$ \Box

\end{document}